\documentclass[preprint,12pt]{elsarticle}

\usepackage{subfigure}
\usepackage{mathrsfs}
\usepackage{eucal}
\usepackage{amsmath}
\usepackage{amssymb}
\usepackage{amsfonts}
\usepackage{graphics,graphicx}

\usepackage{multirow}
\usepackage{lscape}
\usepackage{setspace}

\newcommand{\Eref}[1]{Equation (\ref{#1})}
\newcommand{\fref}[1]{Figure (\ref{#1})}
\newcommand{\Erefs}[1]{Equations (\ref{#1})}
\newcommand{\frefs}[1]{Figures~(\ref{#1})}

\newcommand{\rmd}{\mathrm{d}}

\newcommand{\xx}{\mathbf{x}}

\newcommand{\KK}{\mathbf{K}}
\newcommand{\BB}{\mathbf{B}}
\newcommand{\DD}{\mathbf{D}}

\newcommand{\bvsig}{\boldsymbol{\sigma}}
\newcommand{\bveps}{\boldsymbol{\varepsilon}}
\newcommand{\bn}{\mathbf{N}}
\newcommand{\bfm}{\mathbf{M}}

\journal{Physica E}

\begin{document}

\begin{frontmatter}

\title{On the application of polynomial and NURBS functions for nonlocal response of low dimensional structures}

\author[a]{S~Natarajan\fnref{label2}}
\author[b]{T~Rabczuk}
\author[c]{S~Bordas}
\author[d]{D~Roy Mahapatra}

\address[a]{Post-Doctoral Researcher, Department of Aerospace Engineering, Indian Institute of Science, Bangalore, India.}
\address[b]{Professor, Department of Civil Engineering, Bauhaus-Universit\"{a}t Weimar, Germany}
\address[c]{Professor, Institute of Mechanics and Advanced Materials, Theoretical and Computational Mechanics, Cardiff University, U.K.}
\address[d]{Assistant Professor, Department of Aerospace Engineering, Indian Institute of Science, Bangalore, India.}

\fntext[1]{Department of Aerospace Engineering, Indian Institute of Science, Bangalore - 560012, India. E-mail: sundararajan.natarajan@gmail.com}

\date{\today}

\begin{abstract}
In this paper, the axial vibration of cracked beams, the free flexural vibrations of nanobeams and plates based on Timoshenko beam theory and first-order shear deformable plate theory, respectively, using Eringen's nonlocal elasticity theory is numerically studied. The field variable is approximated by Lagrange polynomials and non-uniform rational B-splines. The influence of the nonlocal parameter, the beam and the plate aspect ratio and the boundary conditions on the natural frequency is numerically studied. The influence of a crack on axial vibration is also studied. The results obtained from this study are found to be in good agreement with those reported in the literature.
\end{abstract}


\end{frontmatter}

\section{Introduction}

The inherent assumption in local elasticity is that dimensions of engineering structure are so much larger than the characteristic dimensions of the microstructure. Thus, the classical continuum theories lack the capability of representing the size effects since they do not include any internal length scale~\citep{marangantisharma2007}. Consequently, these theories fail when the specimen size or the wavelength become comparable with the internal length scales of the material. Several modifications of the classical elasticity formulation have been proposed, such as the strain gradient theory~\citep{aifantis1984}, modified coupled stress theory~\citep{mindlin1964,magao2008} and nonlocal elasticity theory~\citep{eringen1983,eringen1972,eringenedelen1972}. A common feature of these theories is that they include one or several intrinsic length scales. The predictions of these theories reduce to those of local continuum theories where the specimen size is much larger than the internal length scale. The key idea of the nonlocal elastic approaches is that within a nonlocal elastic medium, the particles influence one another not simply by contact forces and heat diffusion but also by long range cohesive forces. In this way, the internal length scale, can be considered in the constitutive equations simply as a material parameter, called nonlocal parameter. One approach to estimate these nonlocal material parameters is by matching the phonon dispersion relation computed by these theories with the lattice dynamics dispersion relation~\citep{chenlee2003,peerlingsfleck2004,marangantisharma2007a}.

It is seen from the literature that the amount of work carried on the application of nonlocal and/or gradient elasticity theory to study the dynamic response of nanobeams is considerable~\citep{peddiesonbuchanan2003,reddy2007,murmupradhan2009,phadikarpradhan2010,roqueferreira2011}. Eringen's nonlocal elasticity theory has been applied by several authors to study axial vibrations~\citep{aydogdu2009} and free transverse vibrations of nanostructures~\citep{reddy2007,aydogdu2009a,phadikarpradhan2010}. Reddy and Pang~\citep{reddypang2008} derived governing equations of motion for Euler-Bernoulli beams and Timoshenko beams using the nonlocal differential relations of Eringen~\citep{eringen1983} and presented closed form solutions for beam static bending, vibration and buckling response with various boundary conditions. Recently, Reddy~\citep{reddy2010}, based on von K\'arm\'an nonlinear strains and Eringen's nonlocal theory, derived governing equations of equilibrium for beams. Phadikar and Pradhan~\citep{phadikarpradhan2010} developed a variational formulation for Euler-Bernoulli beam and Kirchoff plate and employed Hermite and Lagrange polynomials to study the response of beams and plates within nonlocal elasticity framework. Roque \textit{et al.,}~\citep{roqueferreira2011} studied bending, buckling and vibration of Timoshenko nanobeams using collocation techniques. Civalek and Demi~\citep{civalekdemir2011} applied nonlocal Euler-Bernoulli beam theory to study static bending and vibration of microtubules using differential quadrature method. Danesh \textit{et al.,}~\cite{daneshfarajpour2012} employed differential quadrature method to study axial vibration of nanorod.
The nonlocal linear elasticity of Eringen has been applied to study free in-plane vibration and flexural vibration of plates~\citep{murmupradhan2009a,ansarirajabiehfard2010,malekzadehsetoodeh2011}. Reddy~\citep{reddy2010} and Aghabahaei and Reddy~\citep{aghababaeireddy2009} extended classical, first and third order shear deformation theory of plates to account for nonlocal effects. Aghababaei and Reddy~\citep{aghababaeireddy2009} reformulated the third-order shear deformation theory using nonlocal linear theory of elasticity and presented analytical solutions of bending and vibration of plates. The free vibration of single-layered and multi-layered graphene sheets have been studied using first order shear deformation nonlocal plate theory~\citep{ansarirajabiehfard2010,ansarisahmani2010}.

Nanostructures, like any other structures, may develop a flaw during manufacturing. For example, thermally induced cracks are observed in ZnO nanorods during fabrication~\cite{kirkhamwang2008}.The presence of a crack affects the dynamic response of a structural member. This is because, the presence of the crack introduces local flexibility. Moreover, with diminishing external characteristic length, the nonlocal effects cannot be ignored. Recently, Loya \textit{et al.,}~\cite{loyal'opez-puente2009} and Hasheminejad \textit{et al.,}~\cite{hasheminejadgheslaghi2011} studied transverse vibrations of nanobeams based on Euler-Bernoulli theory with and without surface effects, respectively. Hsu \textit{et al.,}~\cite{hsulee2011} studied axial vibration of cracked nanobeam using nonlocal elasticity theory. 

To author's knowledge, the existing numerical approaches are limited to using radial basis functions~\cite{roqueferreira2011} and differential quadrature method~\cite{civalekdemir2011} to study the response of nanostructures based on nonlocal elasticity theory. In this paper, we use Lagrange basis functions and non-uniform rational B splines to approximate the field variable and the response is studied within the framework of the finite element method. The influence of the nonlocal parameter, the boundary conditions and the beam/plate aspect ratio on the fundamental frequency is numerically studied. The free axial vibrations of cracked nanobeams is studied by applying the extended finite element method (XFEM), where the crack is modelled independent of the mesh. This is done by adding suitable 'enrichment' functions to the already existing finite element basis functions. The influence of the boundary conditions, the crack location, the crack severity and the nonlocal parameter on the longitudinal frequency is studied.

The paper is organized as follows. A brief introduction to Eringen's nonlocal elasticity is given in the next section, followed by the governing equations of motion for nanorod, nanobeam based on Timoshenko beam theory and nanoplate based on the first-order shear deformable plate theory (FSDT). Section~\ref{basisfn} presents an overview of the basis functions used in this study and the corresponding discretized form. Section~\ref{numeexamples} presents results for the free axial vibrations of cracked beam and free flexural vibration of Timoshenko beams and first-order shear deformable plates based on nonlocal elasticity theory. The effect of nonlocal parameter, the beam/plate aspect ratio and the boundary conditions on the fundamental frequency is numerically studied, following by concluding remarks in the last section.

\section{Nonlocal elasticity theory}
\label{nonlocal}

Eringen~\citep{eringen1972,eringenedelen1972}, by including the long range cohesive forces, proposed a nonlocal theory in which the stress at a point depends on the strains of an extended region around that point in the body. Thus, the nonlocal stress tensor $\bvsig$ at a point $\xx$ is expressed as:

\begin{equation}
\bvsig = \int_\Omega \alpha(|\xx^\prime - \xx|,\tau) \mathbf{t}(\xx^\prime)~\rmd \xx^\prime 
\label{eqn:eringenInte}
\end{equation}

where $\mathbf{t}(\xx)$ is the classical, macroscopic stress tensor at a point $\xx$ and the kernel function $\alpha(|\xx^\prime - \xx|)$ is called the nonlocal modulus, also referred to as the attenuation or influence function. The nonlocality effects at a reference point $\xx$ produced by the strains at $\xx$ and $\xx^\prime$ are included in the constitutive law by this function. The kernel weights the effects of the surrounding stress states. From the structure of the nonlocal constitutive equation, given by \Eref{eqn:eringenInte}, it can be seen that the nonlocal modulus has the dimension of (length)$^{-3}$~\citep{eringen1972,eringenedelen1972}. Typically, the kernel is a function of the Euclidean distance between the points $\xx$ and $\xx^\prime$. The kernel in \Eref{eqn:eringenInte} has the following properties:

\begin{itemize}
\item The kernel is maximum at $\xx^\prime = \xx$ and decays with $|\xx^\prime -  \xx|$.
\item The classical elasticity limit is included in the limit of vanishing internal characteristic length
\begin{equation}
\lim_{\substack{\tau \rightarrow 0}} \alpha(|\xx^\prime - \xx|) = \delta(|\xx^\prime -  \xx|)
\end{equation}
\item For small internal characteristic length, the nonlocal theory approximates the atomic lattice dynamics.
\item The kernel is a Green's function of a linear differential operator, $L$
\begin{equation}
L \alpha(|\xx^\prime -  \xx|) = \delta(|\xx^\prime -  \xx|)
\end{equation}
\end{itemize}

By appropriate choice of the kernel function, Eringen~\citep{eringen1983} showed that the nonlocal constitute equation given in integral form (\Eref{eqn:eringenInte}) can be represented in an equivalent differential form as:

\begin{equation}
(1 - \tau^2 L^2 \nabla^2) \bvsig = \mathbf{t}
\label{eqn:eringendiffe}
\end{equation}

where $\tau^2 = \frac{\mu}{L^2} = \left( \frac{e_oa}{L} \right)^2$, $e_o$ is a material constant and $a$ and $L$ are the internal and external characteristic lengths, respectively. 

\subsection{Rod}

For nanorod, the constitutive relation is given by~\cite{daneshfarajpour2012}:
\begin{equation}
N - \mu \frac{\partial^2 N}{\partial x^2} = EA \frac{\partial u}{\partial x}
\end{equation}

where $EA$ is the effective axial rigidity and $u$ is the axial displacement. The equations of motion for the axial vibration of nanorod is given by:

\begin{equation}
\frac{\partial}{\partial x} \left( EA \frac{\partial u}{\partial x} \right) = \left(1 - \mu \frac{\partial^2}{\partial x^2} \right) \rho \frac{\partial^2 u}{\partial t^2}
\end{equation}

\subsection{Timoshenko beam theory}

For Timoshenko beam, the constitutive relations are given by~\citep{reddy2010}:

\begin{eqnarray}
M - \mu \frac{\partial^2 M}{\partial x^2} = EI \frac{\partial \phi}{\partial x} \nonumber \\
Q - \mu \frac{\partial^2 Q}{\partial x^2} = GA \kappa_s \left( \phi + \frac{\partial w}{\partial x} \right)
\end{eqnarray}

where $Q$ is the shear force, $G$ is the shear modulus and $\kappa_s$ is the shear correction factor. The displacement field based on the Timoshenko beam theory is given by:

\begin{equation}
u(x,z,t) = u_o(x,t) + z \phi(x,t); \hspace{0.5cm} w(x,z,t) = w_o(x,t).
\end{equation}

where $u_o$ and $w_o$ are the axial and the transverse displacements of the point on the mid-plane (i.e. $z =$ 0) of the beam. The nonzero strain of the Timoshenko beam theory is

\begin{equation}
\varepsilon_{xx} = \frac{ \partial u_o}{\partial x} + z \frac{\partial \phi}{\partial x}; \hspace{0.5cm} \gamma_{xz} = \frac{ \partial w_o}{\partial x} + \phi.
\end{equation}

The Timoshenko beam theory requires shear correction factors to compensate for the error due to the constant shear stress assumption. The governing equations of motion in the case of nonlocal elasticity is given by:

\begin{eqnarray}
\frac{\partial}{\partial x} \left[ G A \kappa_s \left( \phi + \frac{\partial w}{\partial x} \right) \right] = m_o \left( \frac{\partial ^2 w}{\partial t^2} - \mu \frac{\partial^4 w}{\partial x^2 \partial t^2} \right) \nonumber \\
\frac{\partial}{\partial x} \left( EI \frac{\partial \phi}{\partial x}\right) - GA\kappa_s \left( \phi + \frac{\partial w}{\partial x} \right) = m_2 \left( \frac{\partial ^2 \phi}{\partial t^2} - \mu  \frac{\partial^4 \phi}{\partial x^2 \partial t^2} \right) 
\label{eqn:tbt_governeqn}
\end{eqnarray}

\subsection{Reissner-Mindlin plate}
Consider a rectangular Cartesian coordinate system $(x,y,z)$ with $xy-$ plane coinciding with the undeformed middle plane of the plate and the $z$- coordinate is taken positive downwards. The displacement field based on first order shear deformation plate theory (FSDT), also referred to as the Mindlin plate theory is given by:

\begin{eqnarray}
u(x,y,z,t) &=& u_o(x,y,t) + z \theta_x(x,y,t) \nonumber \\
v(x,y,z,t) &=& v_o(x,y,t) + z \theta_y(x,y,t) \nonumber \\
w(x,y,z,t) &=& w_o(x,y,t) 
\label{eqn:displacements}
\end{eqnarray}

\begin{equation}
\bveps  = \left\{ \begin{array}{c} \bveps_p \\ 0 \end{array} \right \}  + \left\{ \begin{array}{c} z \bveps_b \\ \bveps_s \end{array} \right\} 
\label{eqn:strain1}
\end{equation}

The midplane strains $\bveps_p$, bending strain $\bveps_b$, shear strain $\bveps_s$ in \Eref{eqn:strain1} are written as

\begin{equation}
\renewcommand{\arraystretch}{1.5}
\bveps_p = \left\{ \begin{array}{c} u_{o,x} \\ v_{o,y} \\ u_{o,y}+v_{o,x} \end{array} \right\}, \hspace{0.5cm}
\renewcommand{\arraystretch}{1.5}
\bveps_b = \left\{ \begin{array}{c} \theta_{x,x} \\ \theta_{y,y} \\ \theta_{x,y}+\theta_{y,x} \end{array} \right\} \hspace{0.5cm}
\renewcommand{\arraystretch}{1.5}
\bveps_s = \left\{ \begin{array}{c} \theta _x + w_{o,x} \\ \theta _y + w_{o,y} \end{array} \right\}.
\renewcommand{\arraystretch}{1.5}
\end{equation}

where the subscript `comma' represents the partial derivative with respect to the spatial coordinate succeeding it. The equations of motion is given by:

\begin{eqnarray}
\frac{\partial N_{xx}}{\partial x} + \frac{\partial N_{xy}}{\partial x} = I_{11} \ddot{u} + I_{12} \ddot{\theta_x} \nonumber \\
\frac{\partial N_{xy}}{\partial x} + \frac{\partial N_{yy}}{\partial x} = I_{11} \ddot{v} + I_{12} \ddot{\theta_y} \nonumber \\
\frac{\partial N_{xz}}{\partial x} + \frac{\partial N_{yz}}{\partial x} = I_{11} \ddot{w} \nonumber \\
\frac{\partial M_{xx}}{\partial x} + \frac{\partial M_{xy}}{\partial x} = I_{12} \ddot{u} + I_{22} \ddot{\theta_x} \nonumber \\
\frac{\partial M_{xy}}{\partial x} + \frac{\partial M_{yy}}{\partial x} = I_{12} \ddot{v} + I_{22} \ddot{\theta_y}
\label{eqn:plateeqeqn}
\end{eqnarray}

where $(I_{11}, I_{12}, I_{22}) = \int\limits_{-h/2}^{h/2} \rho(1,z,z^2)~\rmd z$. The nonlocal membrane stress resultants $\bn$ and the bending stress resultants $\bfm$ can be related to the membrane strains, $\bveps_p$ and bending strains $\bveps_b$ through the following constitutive relations~\citep{reddy2010}:

\begin{eqnarray}
\bn &=& \left\{ \begin{array}{c} N_{xx} \\ N_{yy} \\ N_{xy} \end{array} \right\}(1- \mu \nabla^2) = \mathbf{A}_e \bveps_p + \BB_{be} \bveps_b \nonumber \\
\bfm &=& \left\{ \begin{array}{c} M_{xx} \\ M_{yy} \\ M_{xy} \end{array} \right\}(1- \mu \nabla^2) = \BB_{be} \bveps_p + \DD_b \bveps_b  
\label{eqn:platenonlocal}
\end{eqnarray}

where the matrices $\mathbf{A}_e = A_{ij}, \BB_{be} = B_{ij}$ and $\DD_b = D_{ij}; (i,j=1,2,6)$ are the extensional, bending-extensional coupling and bending stiffness coefficients and are defined as

\begin{equation}
\left\{ A_{ij}, ~B_{ij}, ~ D_{ij} \right\} = \int_{-h/2}^{h/2} \overline{Q}_{ij} \left\{1,~z,~z^2 \right\}~dz, \hspace{0.2cm} i=1,2,6
\end{equation}

Similarly, the transverse shear force $Q = \{Q_{xz},Q_{yz}\}$ is related to the transverse shear strains $\bveps_s$ through the following equation

\begin{equation}
Q_{ij} = E_{ij} \bveps_s
\end{equation}

where $E_{ij} = \mathbf{E} = \int_{-h/2}^{h/2} \overline{Q} \upsilon_i \upsilon_j~dz;~ (i,j=4,5)$ is the transverse shear stiffness coefficient, $\upsilon_i, \upsilon_j$ is the transverse shear correction factors for non-uniform shear strain distribution through the plate thickness. The stiffness coefficients $\overline{Q}_{ij}$ are defined as

\begin{eqnarray}
\overline{Q}_{11} = \overline{Q}_{22} = {E(z) \over 1-\nu^2}; \hspace{1cm} \overline{Q}_{12} = {\nu E(z) \over 1-\nu^2}; \hspace{1cm} \overline{Q}_{16} = \overline{Q}_{26} = 0 \nonumber \\
\overline{Q}_{44} = \overline{Q}_{55} = \overline{Q}_{66} = {E(z) \over 2(1+\nu) }
\end{eqnarray}

The equations of motion in terms of midplane displacements and rotations are obtained by substituting the nonlocal constitutive relations (\Eref{eqn:platenonlocal}) into the equations of motion~\Eref{eqn:plateeqeqn}. For detailed derivation, interested readers are referred to the literature~\citep{murmupradhan2009a,reddy2010,malekzadehsetoodeh2011}.



\section{Basis functions and spatial discretization}\label{basisfn}

In this study, the finite element model has been developed using Lagrange polynomials, NURBS basis function. In this section, a brief overview of different basis functions used in this study is given. As it is not the scope of this study to review the different basis functions, interested readers are referred to the literature and references therein.

\paragraph{Lagrange polynomials}
The Lagrange polynomials~\citep{zienkiewicztaylor2000} for a two noded element is given by:

\begin{equation}
N_1 = 1 - \frac{x}{\ell}; \hspace{0.5cm} N_2 = \frac{x}{\ell}
\end{equation}

where $\ell$ is the length of the element. The Lagrange polynomials are used to model the nanorod.

\paragraph{Field consistent Q4 and Q8 element}
The plate element employed here is a $\mathcal{C}^0$ continuous shear flexible field consistent element with five degrees of freedom $(u_o,v_o,w_o,\theta_x,\theta_y)$ at four nodes in an 4-noded quadrilateral (Q4) element and eight nodes in a 8-noded quadrilateral element (Q8). If the interpolation functions for Q4/Q8 are used directly to interpolate the five variables $(u_o,v_o,w_o,\theta_x,\theta_y)$ in deriving the shear strains and membrane strains, the element will lock and show oscillations in the shear and membrane stresses. The field consistency requires that the transverse shear strains and membrane strains must be interpolated in a consistent manner. Thus, the $\theta_x$ and $\theta_y$ terms in the expressions for shear strain $\bveps_s$ have to be consistent with the derivative of the field functions, $w_{o,x}$ and $w_{o,y}$. This is achieved by using field redistributed substitute shape functions to interpolate those specific terms, which must be consistent as described in~\cite{somashekarprathap1987,ganapathivaradan1991}. This element is free from locking and has good convergence properties. For complete description of the element, interested readers are referred to the literature~\cite{somashekarprathap1987,ganapathivaradan1991}, where the element behaviour is discussed in great detail. Since the element is based on the field consistency approach, exact integration is applied for calculating various strain energy terms.

\paragraph{Non-uniform Rational B-Splines}
The key ingredients in the construction of NURBS basis functions are: the knot vector (a non decreasing sequence of parameter value, $\xi_i \le \xi_{i+1}, i = 0,1,\cdots,m-1$), the control points, $P_i$, the degree of the curve $p$ and the weight associated to a control point, $w$. The i$^{th}$ B-spline basis function of degree $p$, denoted by $N_{i,p}$ is defined as~\cite{piegltiller1996,hughescottrell2005}:

\begin{eqnarray}
N_{i,0}(\xi) = \left\{ \begin{array}{cc} 1 & \textup{if} \hspace{0.2cm} \xi_i \le \xi \le \xi_{i+1} \\
0 & \textup{else} \end{array} \right. \nonumber \\
N_{i,p}(\xi) = \frac{ \xi- \xi_i}{\xi_{i+p} - \xi_i} N_{i,p-1}(\xi) + \frac{\xi_{i+p+1} - \xi}{\xi_{i+p+1}-\xi_{i+1}}N_{i+1,p-1}(\xi)
\end{eqnarray}

A $p^{th}$ degree NURBS curve is defined as follows:

\begin{equation}
\mathbf{C}(\xi) = \frac{\sum\limits_{i=0}^m N_{i,p}(\xi)w_i \mathbf{P}_i}{\sum\limits_{i=0}^m N_{i,p}(\xi)w_i}
\end{equation}

where $\mathbf{P}_i$ are the control points and $w_i$ are the associated weights.

The displacement field variable is approximated by one of the above mentioned basis functions. By multiplying the governing equations of motion by appropriate test functions and by employing Bubnov-Galerkin method, we get the following set of algebraic equations:

\begin{equation}
\bfm \ddot{\boldsymbol{\delta}} + \KK \boldsymbol{\delta} = 0
\end{equation}

where $\bfm$ is the consistent mass matrix, $\KK$ is the stiffness matrix and $\boldsymbol{\delta}$ is the vector of the degree of freedom associated to the displacement field. After substituting the characteristic of the time function $\ddot{\boldsymbol{\delta}} = \omega^2 \boldsymbol{\delta}$, the following algebraic equation is obtained:

\begin{equation}
( \KK - \omega^2 \bfm) \boldsymbol{\delta}  = \mathbf{0}
\end{equation}

where $\omega$ is the natural frequency. The frequencies are obtained using the standard generalized eigenvalue algorithm.

\section{Numerical examples}
\label{numeexamples}
In this section, we present the axial vibrations of cracked beam and natural frequencies for Timoshenko beam and first-order shear deformable plate, based on Eringen's nonlocal elasticity theory. In all cases, we present the non dimensionalized free flexural frequencies as, unless specified otherwise:

\paragraph*{For axial vibration of beams}
\begin{equation}
\Omega = \omega  L \sqrt{ \frac{\rho}{EA}}
\end{equation}

where $L, \rho$ and $EA$ denote the length of the beam, mass density and effective axial rigidity.

\paragraph*{For flexural vibration of Timoshenko beam}
\begin{equation}
\Omega = \omega L^2 \sqrt{ \frac{\rho A}{EI}}
\end{equation}

where $EI$ is the flexural rigidity and $\rho$ is the mass density. 

\paragraph*{For plates}
\begin{equation}
\Omega = \omega h \sqrt{ \frac{\rho}{G}}
\end{equation}

where $\omega$ is the natural frequency. In all cases, the effect of the nonlocal parameter is studied by the frequency ratio, which is defined as:

\begin{equation}
\textup{Frequency ratio} = \frac{\Omega_{NL}}{\Omega_L}
\end{equation}

The frequency ratio is a measure of the error made by neglecting the nonlocal effects. For example, if $\frac{\Omega_{NL}}{\Omega_L}=$ 1, the effect of nonlocal parameter is not significant and becomes important for any other value.

\subsection{Axial vibration of cracked beams}
Consider a slender nanorod with uniform cross section $A$ along the length $L$. The crack at distance of $C$ from the left end is simulated by an equivalent linear spring. The analytical expression for the eigen frequencies are given by~\cite{hsulee2011}:

\paragraph*{Clamped-free} 

\begin{equation}
\cos \left( \frac{\beta}{\sqrt{1-\mu \beta^2}} \right) - K \beta \sqrt{1-\mu \beta^2} \cos \left( \frac{b \beta}{\sqrt{1-\mu \beta^2}} \right) \sin \left( \frac{\beta - b\beta}{\sqrt{1-\mu \beta^2}} \right) = 0
\label{eqn:clamfreeana}
\end{equation}

\paragraph*{Clamped-Clamped}

\begin{equation}
2\sin \left( \frac{\beta}{\sqrt{1 - \mu \beta^2}} \right) + K \beta \sqrt{1 - \mu \beta^2} \left[ \cos \left( \frac{\beta}{\sqrt{1 - \mu \beta^2}} \right) + \cos \left( \frac{b - 2b \beta}{\sqrt{1 - \mu \beta^2}} \right) \right] = 0
\label{eqn:clamclamana}
\end{equation}

when $K (= \frac{EA}{kL}), b$ and $\mu (= \left(\frac{e_oa}{L}\right)^2 )$ are known, the dimensionless longitudinal frequency $\Omega$ of a single cracked nanobeam with nonlocal effect is solved by using \Erefs{eqn:clamfreeana} - (\ref{eqn:clamclamana}). For a nanobeam without a crack, i.e., $K=$ 0, the various vibration modes are expressed as~\cite{aydogdu2009,filizaydogdu2010,hsulee2011}:

\paragraph*{Clamped-free}

\begin{equation}
\beta_{\rm{NL}} = \frac{ (2n-1) \frac{\pi}{2}} {\sqrt{ 1 + (2n-1)^2 \mu \left( \frac{\pi}{2} \right)^2}}, \hspace{0.2cm} n = 1,2,3,\cdots
\label{eqn:nonlocalclamfree}
\end{equation}

\paragraph*{Clamped-clamped}

\begin{equation}
\beta_{\rm{NL}} = \frac{ n \pi }{\sqrt{ 1 + n^2 \mu \left( \pi \right)^2}}, \hspace{0.2cm} n = 1,2,3,\cdots
\label{eqn:nonlocalclamclamp}
\end{equation}

\begin{figure}[htpb]
\centering
\subfigure[$e_oa/L = $ 0]{\includegraphics[scale=0.7]{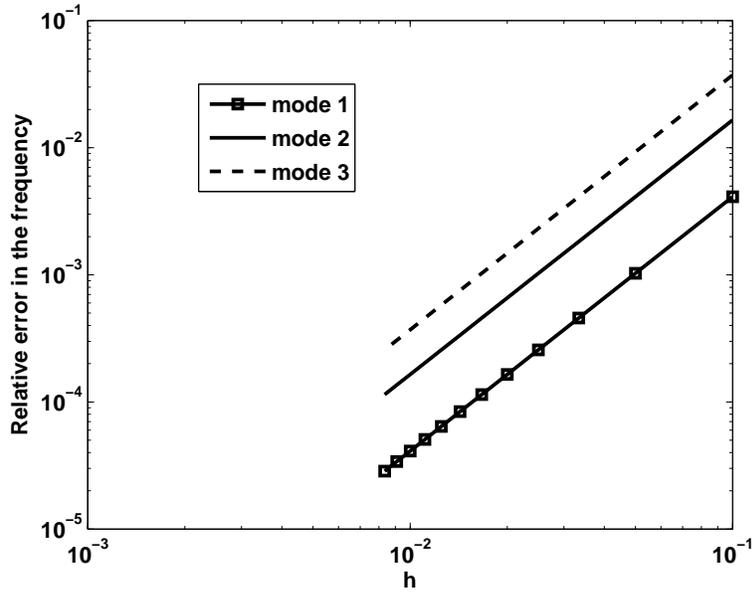}}
\subfigure[$e_oa/L = $ 0.1]{\includegraphics[scale=0.7]{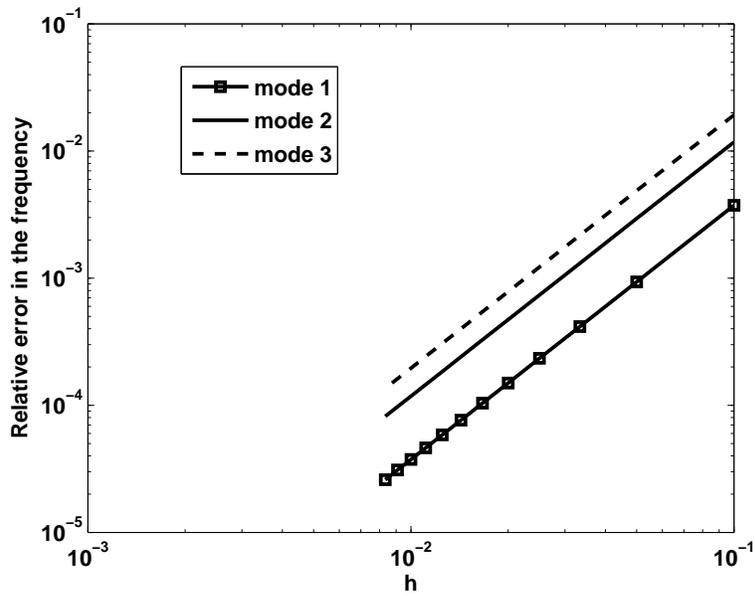}}
\caption{Convergence of the first three modes with mesh size.}
\label{fig:meshconve}
\end{figure}

\begin{figure}[htpb]
\centering
\subfigure[clamped-free]{\includegraphics[scale=0.7]{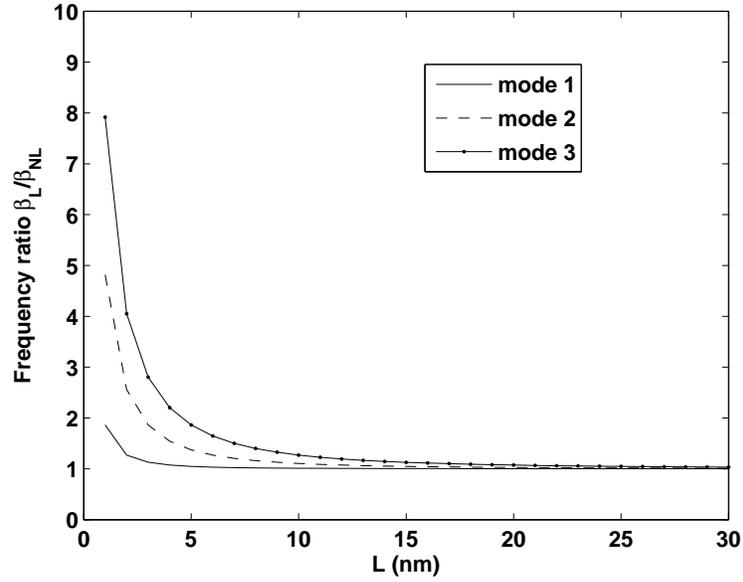}}
\subfigure[clamped-clamped]{\includegraphics[scale=0.7]{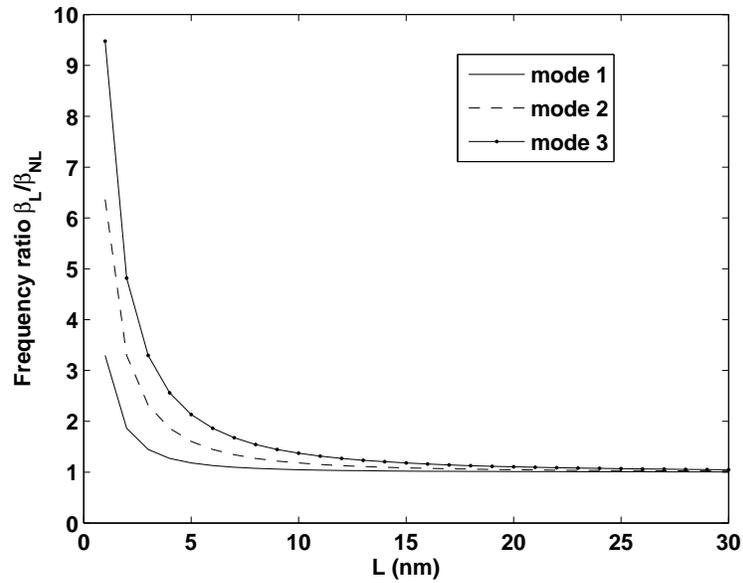}}
\caption{Frequency ratio $(\beta_{\rm L}/\beta_{\rm {nl}})$ as a function of beam length with the nonlocal parameter $e_oa =$ 1nm. It can be seen that as the length of the beam increases with respect to the characteristic internal length, the nonlocal effects decreases and for very large beams, the nonlocal effects are negligible.}
\label{fig:effectofLength}
\end{figure}

By choosing $\mu=$ 0 in \Erefs{eqn:nonlocalclamfree} and (\ref{eqn:nonlocalclamclamp}), the various vibration modes for local elasticity are obtained. \fref{fig:meshconve} shows the convergence of first three modes with mesh size with and without nonlocal effects. It can be seen that with decreasing mesh size, the analytical solution is approached. Based on a progressive mesh refinement, 100 elements is found to be adequate for the study. Tables \ref{clampedFree_Conve} - \ref{clampedclamped_comp} gives a comparison of computed frequencies for clamped-free and clamped-clamped nanobeam with a single crack, respectively. It can be seen that the numerical results from the present formulation are found to be in good agreement with the existing solutions. The effect of nonlocal parameter on the fundamental frequencies is shown in \fref{fig:effectofLength} for a fixed nonlocal parameter $e_oa = $1 nm. In this case, the length of the beam is varied. It can be seen that as the length of the beam increases with respect to the characteristic internal length, the nonlocal effects decreases and for very large beams, the nonlocal effects are negligible.

\begin{table}[htpb]
\centering
\renewcommand\arraystretch{1.5}
\caption{Comparison of natural frequencies of the clamped-free beam with a single crack located at $C/L=$ 0.2002, crack severity $K=$ 0.1144 and nonlocal parameter $e_oa/L=$ 0.}
\begin{tabular}{clll}
\hline
Mode & XFEM & Ref.~\cite{hsulee2011} & Ref.~\cite{singh2009}\\
\hline
1 & 1.4228 & 1.4278 & 1.4278 \\
2 & 4.4429 & 4.5579 & 4.5576 \\
3 & 7.8559 & 7.8540 & 8.8540 \\
4 & 10.4289 & 10.4471 & 10.4486 \\
\hline
\end{tabular}
\label{clampedFree_Conve}
\end{table}

\begin{table}[htpb]
\centering
\renewcommand\arraystretch{1.5}
\caption{Comparison of natural frequencies of the clamped-clamped beam with a single crack located at $C/L=$ 0.25 for various internal length and crack parameter~\cite{hsulee2011}.}
\begin{tabular}{crrrrrr}
\hline
$e_oa/L$ & \multicolumn{2}{c}{$K=$ 0.065} & \multicolumn{2}{c}{$K=$ 0.35} & \multicolumn{2}{c}{$K=$ 2} \\
\cline{2-7}
& XFEM & Ref.~\cite{hsulee2011} & XFEM & Ref.~\cite{hsulee2011} & XFEM & Ref.~\cite{hsulee2011} \\
\hline
0.2 & 2.6144 & 2.6173 & 2.4649 & 2.4668 & 2.1503 & 2.1506 \\
0.4 & 1.9455 & 1.9467 & 1.9060 & 1.9071 & 1.7660 & 1.7663 \\
\hline
\end{tabular}
\label{clampedclamped_comp}
\end{table}

\begin{figure}[htpb]
\centering
\subfigure[clamped-free]{\includegraphics[scale=0.7]{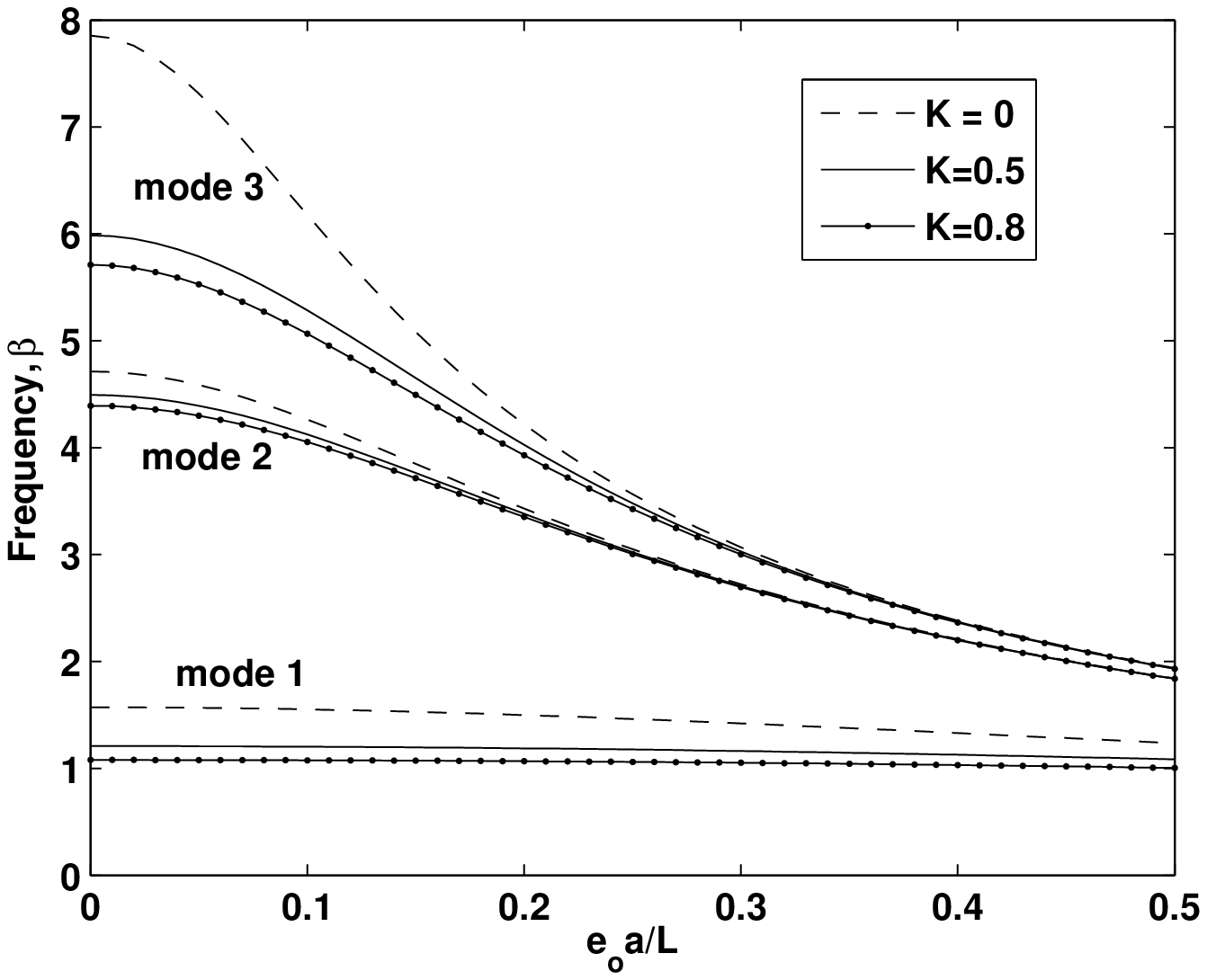}}
\subfigure[clamped-clamped]{\includegraphics[scale=0.7]{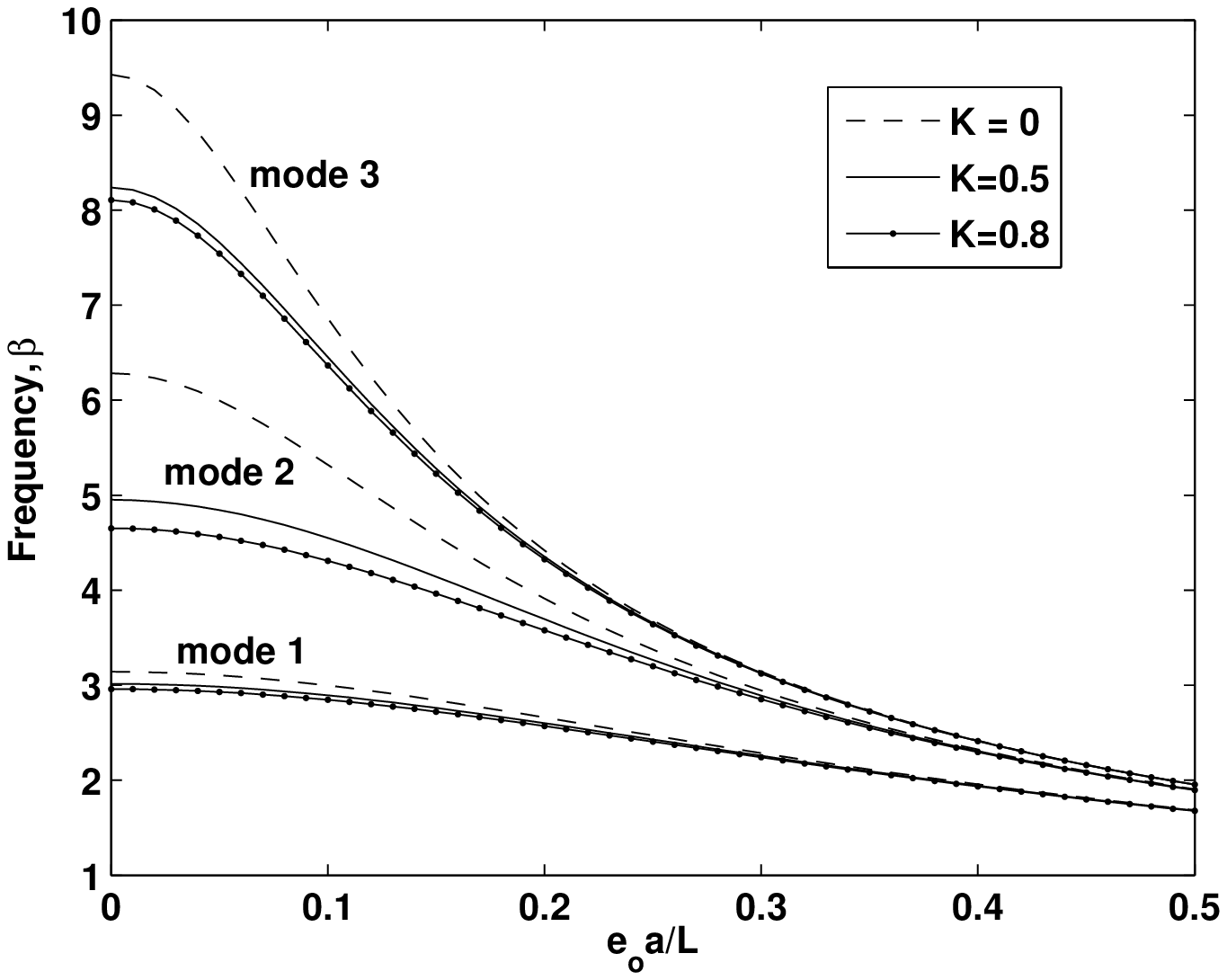}}
\caption{The dimensionless frequency for different crack parameters for the first three vibration modes. The crack is located at a distance 0.4 from the left end.}
\label{fig:effectNL}
\end{figure}

\begin{figure}[htpb]
\centering
\subfigure[]{\includegraphics[scale=0.7]{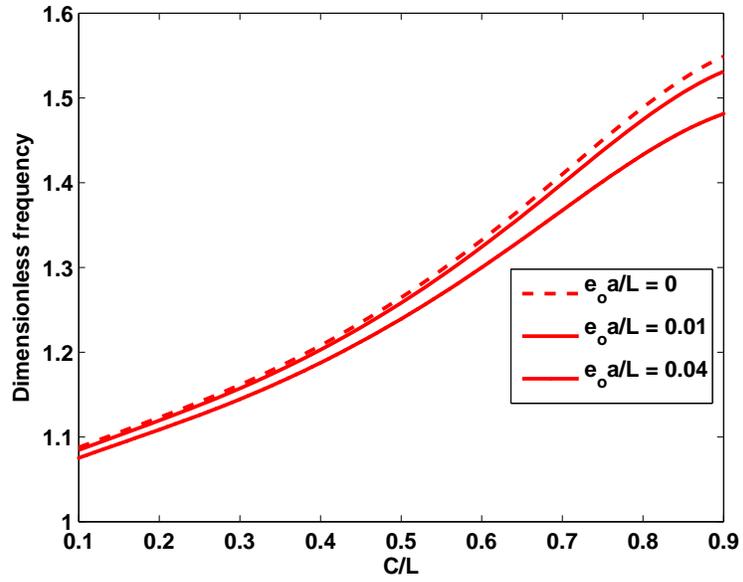}}
\subfigure[]{\includegraphics[scale=0.7]{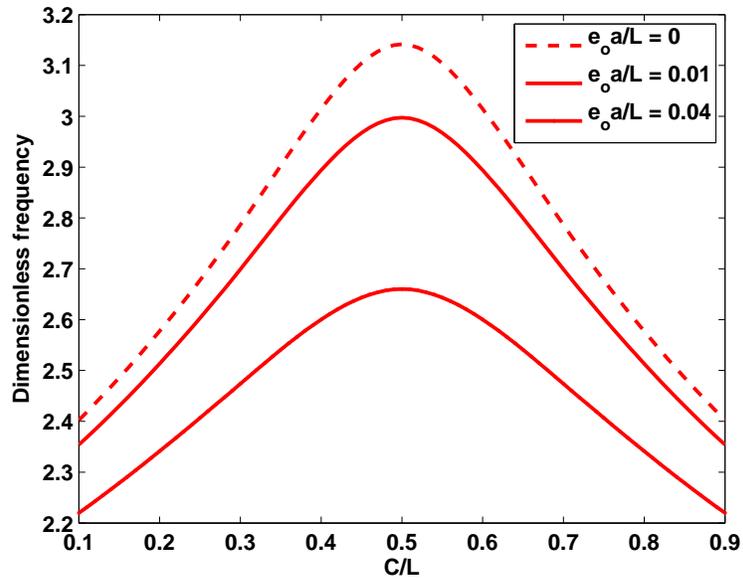}}
\caption{The dimensionless frequency of mode 1 of (a) clamped-free beam and (b) clamped-clamped nanobeam.}
\label{fig:effectCrkloc}
\end{figure}

\fref{fig:effectNL} depicts the influence of boundary conditions and nonlocal parameter on the first three modes of the clamped-free and clamped-clamped beam with a single crack. The crack is located at $C/L=$ 0.4 from the left end. It can be seen that the frequencies are higher for the beam with no crack and with increasing crack parameter $K$ the frequency decreases. The presence of the crack introduces local flexibility and the effect of crack parameter on the frequencies is significant for lower values of the nonlocal parameter $e_oa$. This observation is consistent with the existing results~\cite{hsulee2011}. The effect of crack location along the length of the beam and the boundary conditions on the frequency is shown in \fref{fig:effectCrkloc}. In both cases, the natural frequency of the beam is influenced by the location of the crack. In the case of clamped-free boundary condition, the crack near the free end shows a stronger influence than the crack at the fixed end, while, in the case of clamped-clamped boundary condition, the crack in the middle of the beam has stronger influence. Due to the symmetric boundary conditions, the frequency distribution in the case of clamped-clamped boundary condition is observed. This again is consistent with the results available in the literature~\cite{hsulee2011}.

\subsection{Free flexural vibration of beams}
In the next example, consider a nanobeam of length $a=$ 10 and thickness $h$ and with following material properties: Young's modulus, $E=$ 30 $\times$ 10$^6$, Poisson's ratio $\nu=$ 0.3 and mass density $\rho=$ 1. The value of the shear correction factor $\kappa_s$ is taken to be $5/6$.

\begin{table}[htpb]
\centering
\renewcommand\arraystretch{1.5}
\caption{Comparison of natural frequencies of simply supported beams based on Timoshenko beam theory: $^\diamond$\cite{reddy2007}, $^{\dagger\dagger}$1501 global collocation points~\cite{roqueferreira2011}, $^\flat$NURBS basis functions with $p =$ 3 and 82 elements, $^\ast$Lagrange polynomials (100 elements).}
\begin{tabular}{ccllll}
\hline
$a/h$ & $\mu$ & & \multicolumn{3}{c}{Timoshenko beam}\\
\cline{3-6}
&  & Ref.$^\diamond$ & RBF$^{\dagger\dagger}$ & NURBS$^\flat$ & FEM$^\ast$ \\
\hline
\multirow{3}{*}{100}& 0 &  9.8683 & 9.8283 & 9.8680 & 9.8630 \\
& 1 &  9.4147 & 9.3766 & 9.4144 & 9.4096 \\
& 5 &  8.0750 & 8.0423 & 8.0748 & 8.0706\\
\cline{2-6}
\multirow{3}{*}{20}& 0 & 9.8381 & 9.8059 & 9.8281 & 9.7955\\
& 1 &  9.3858 & 9.3551 & 9.3763 & 9.3452\\
& 5 &  8.0503 & 8.0239 &  8.0421 & 8.0154 \\
\cline{2-6}
\multirow{3}{*}{10}& 0 & 9.7454 & 9.7792 &  9.7075 & 9.5886 \\
& 1 & 9.2973 & 9.3294 & 9.2612 &  9.1460\\
& 5 & 7.9744 & 8.0014 & 7.9434 &  7.8445\\
\hline
\end{tabular}
\label{TBT_comparison}
\end{table}

\begin{figure}[htpb]
\centering
\includegraphics[scale=0.6]{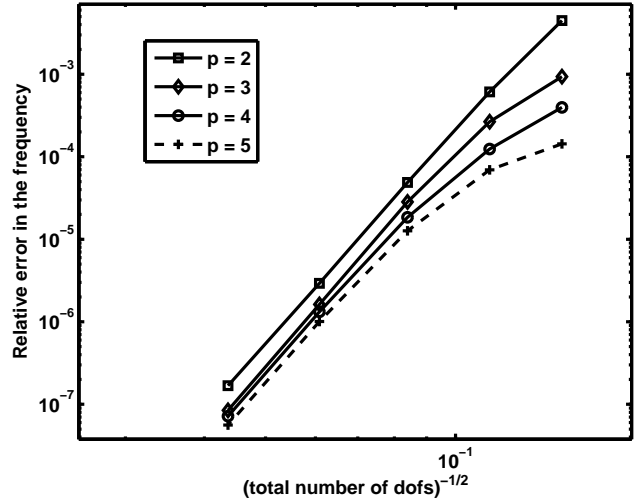}
\caption{Convergence of the first mode with mesh size for various degree of the curve.}
\label{fig:nurbsconvetbt}
\end{figure}
\begin{figure}[htpb]
\centering
\includegraphics[scale=0.6]{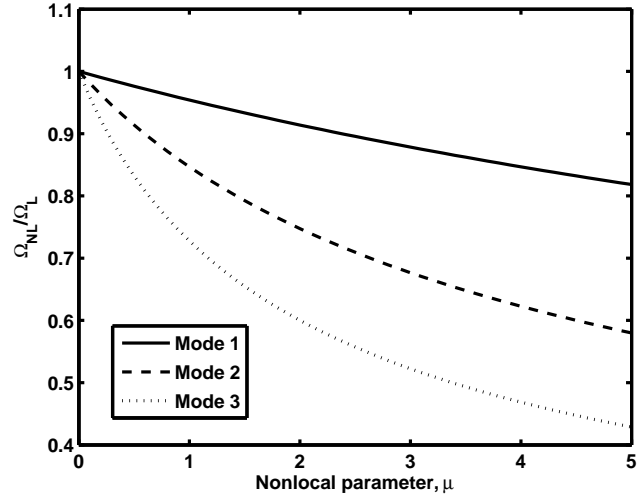}
\caption{Frequency ratio $(\Omega_{NL}/\Omega_L)$ as a function of nonlocal parameter $\mu$ for $a/h=$ 100 and $L=$ 10 for Timoshenko beam. It can be seen that as the nonlocal parameter increases, the nonlocal effects are prominent. For a fixed nonlocal parameter, the effect of nonlocal parameter on higher frequencies is greater than that for lower modes.}
\label{fig:timointlenparam}
\end{figure}
\begin{figure}[htpb]
\centering
\includegraphics[scale=0.6]{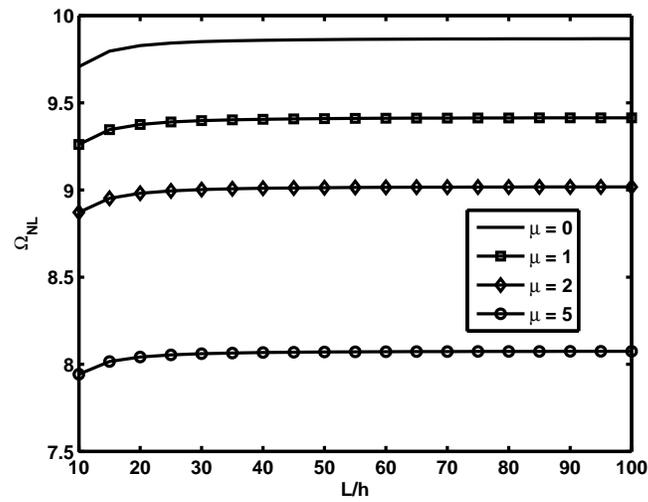}
\caption{Mode 1 non-dimensionalized frequency $(\Omega_{NL})$ as a function of beam aspect ratio $L/h$ for different nonlocal parameters for Timoshenko beam. It can be seen that the effect of nonlocal parameter is to decrease the frequency, whereas, increasing the beam aspect ratio, increases the fundamental frequency}
\label{fig:timoahintlenparam}
\end{figure}

The numerical results for free vibration of Timoshenko beam are given in Table~\ref{TBT_comparison}. It can be seen that the numerical results from the present study are in good agreement with the existing solutions~\citep{reddy2007,roqueferreira2011}. \fref{fig:nurbsconvetbt} shows the convergence in case of the NURBS basis function, for mode 1 frequency with decreasing mesh size for various order of the curve. It can be seen that increasing the order of the curve, yields more accurate results. With decreasing mesh size, all the results approach analytical solution. Note that, the accuracy on a coarse mesh can be improved by increasing the order of the curve. The effect of nonlocal parameter and beam aspect ratio is shown in \frefs{fig:timointlenparam} and (\ref{fig:timoahintlenparam}), respectively. It can be observed that increasing the nonlocal parameter decreases the fundamental frequency for both Euler-Bernoulli beam and Timoshenko beam. The nonlocal parameter has greater influence on higher modes. The variation of frequency ratio $\Omega_{NL}/\Omega_L$ is shown in \fref{fig:timointlenparam}. For mode 1, the variation is almost linear, but for higher modes, the variation becomes nonlinear. It can be seen from \fref{fig:timoahintlenparam} that increasing the beam aspect ratio increases the frequencies, while increasing the nonlocal parameter decreases the frequency. This is because, increasing the aspect ratio, increases the flexibility of the structure, on the other hand, the nonlocal parameter decreases the flexibility. The effect of various boundary condition on the fundamental frequency is illustrated in Table \ref{TBT_variousbc}. It can be seen that the effect of nonlocal parameter is to decrease the frequency. The clamped-clamped boundary condition increases the fundamental frequency, while the clamped-free lowers the frequency compared to the simply supported case.

\begin{table}[htpb]
\centering
\renewcommand\arraystretch{1.5}
\caption{Non-dimensional frequency of Timoshenko beam under various boundary conditions.}
\begin{tabular}{cclll}
\hline
$a/h$ & $\mu$ & \multicolumn{3}{c}{NURBS, $p=$ 3}\\
\cline{3-5}
&  & SS & CC & CF \\
\hline
\multirow{4}{*}{100}& 0 & 9.8680 & 22.3892 & 3.5178\\
& 1 & 9.4144 & 21.1228 & 3.4385 \\
& 2 & 9.0180 & 20.0450 & 3.3639\\
& 5 & 8.0748 & 17.5791 & 3.1646 \\
\cline{2-5}
\multirow{4}{*}{20}& 0 & 9.8281 & 21.9967 & 3.5091\\
& 1 & 9.3763 & 20.7595 & 3.4303 \\
& 2 & 8.9816 & 19.7049 & 3.3561 \\
& 5 & 8.0421 & 17.2877 & 3.1577 \\
\cline{2-5}
\multirow{4}{*}{10}& 0 & 9.7075 & 20.9726 & 3.4884\\
& 1 & 9.2612 & 19.8083 & 3.4107 \\
& 2 & 8.8713 & 18.8124 & 3.3375 \\
& 5 & 7.9434 & 16.5200 & 3.1415 \\
\hline
\end{tabular}
\label{TBT_variousbc}
\end{table}

\subsection{Free vibration of plates}

Consider a plate of uniform thickness, $h$ and with length and width as $a$ and $b$, respectively. In this case, the displacement field is approximated by Lagrange elements (Q4 and Q8) and with NURBS basis functions. The computed frequencies for a square simply supported plate is given in Table~\ref{table:SSisoplate_comp}. It can be seen that the numerical results from the present study are found to be in good agreement with the existing solutions. Table~\ref{table:SSisoplate_comp} gives the results of non-dimensional frequency for first-order nonlocal plate theory for different values of nonlocal parameter, the plate aspect ratio and for various basis functions. Again, it can be seen that the nonlocal theory predicts smaller values of natural frequency than the local theory. 

\begin{table}[htpb]
\renewcommand\arraystretch{1.5}
\caption{Comparison of fundamental frequency $\omega h \sqrt{\rho/G}$ for simply supported plate. $^{\dagger\dagger}$\cite{aghababaeireddy2009}, $^\Box$mesh size 40 $\times 40$, $^\ast$mesh size 8 $\times$ 8, $^\sharp$Order of curve = 3, 5 $\times$ 5 control points.} 
\centering
\begin{tabular}{cccllll}
\hline
$a/b$ & $a/h$ & $\mu$ & Ref.$^{\dagger\dagger}$ & \multicolumn{3}{c}{Method} \\
\cline{5-7}
&  &  &  & Q4$^\Box$ & Q8$^\ast$ & NURBS$^\sharp$\\
\hline
\multirow{6}{*}{1} & \multirow{3}{*}{10} & 0 & 0.0930 & 0.0927 & 0.0926 & 0.0929\\
& & 1 & 0.0850 & 0.0847& 0.0846 &0.0849\\
& & 5 & 0.0660 & 0.0657 & 0.0657 &0.0659\\
\cline{3-7}
& \multirow{3}{*}{20} & 0 & 0.0239 & 0.0240 & 0.0238 &0.0239\\
& & 1 & 0.0218 & 0.0219 & 0.0218 &0.0219\\
& & 5 & 0.0169 & 0.0170 & 0.0169 &0.0170\\
\hline
\multirow{6}{*}{2} & \multirow{3}{*}{10} & 0 & 0.0589 & 0.0588 & 0.0587 &0.0590\\
& & 1 & 0.0556 & 0.0554 & 0.0554 &0.0556\\
& & 5 & 0.0463 & 0.4620 & 0.0462 &0.0464\\
\cline{3-7}
& \multirow{3}{*}{20} & 0 & 0.0150 & 0.0150& 0.0150 & 0.0151\\
& & 1 & 0.0141 & 0.0141& 0.0141 &0.0141\\
& & 5 & 0.0118 & 0.0118& 0.0118 & 0.0118\\
\hline
\end{tabular}
\label{table:SSisoplate_comp}
\end{table}

As another example, the nonlocal frequency -to- local frequency ratio $(\Omega_{N}/\Omega_L)$ is computed for a simply supported isotropic square plate and the results are compared with those available in the literature~\citep{pradhanphadikar2009,malekzadehsetoodeh2011}. Frequency ratio for different values of nonlocal parameter are presented in Table~\ref{table:freqratiovalid}. It can be seen that the results from the present formulation are in good agreement with the existing results.

\begin{table}[htpb]
\renewcommand\arraystretch{1.5}
\caption{Comparison of frequency ratio for a simply supported square plate ($a=$ 10, $h=$ 0.34). $^{\dagger\dagger}$\cite{aghababaeireddy2009}, $^{\circ}$\cite{malekzadehsetoodeh2011}, $^\Box$mesh size 40 $\times 40$, $^\ast$mesh size 8 $\times$ 8, $^\sharp$Order of curve = 3, 5 $\times$ 5 control points.} 
\centering
\begin{tabular}{clllll}
\hline
$\mu$ & Ref.$^{\diamond}$ & Ref.$^{\circ}$ & \multicolumn{3}{c}{Method} \\
\cline{4-6}
 &  &  & Q4$^\Box$ & Q8$^\ast$ & NURBS$^\sharp$ \\
\hline
0 & 1.0000 & 1.0000 & 1.0000 & 1.0000 & 1.0000 \\
1 & 0.9139 & 0.9139 & 0.9099 & 0.9107 & 0.9107 \\
2 & 0.8467 & 0.8467 & 0.8393 & 0.8468 & 0.8468 \\
3 & 0.7925 & 0.7925 & 0.7857 & 0.7928 & 0.7928 \\ 
\hline
\end{tabular}
\label{table:freqratiovalid}
\end{table}

\begin{figure}[htpb]
\centering
\includegraphics[scale=0.6]{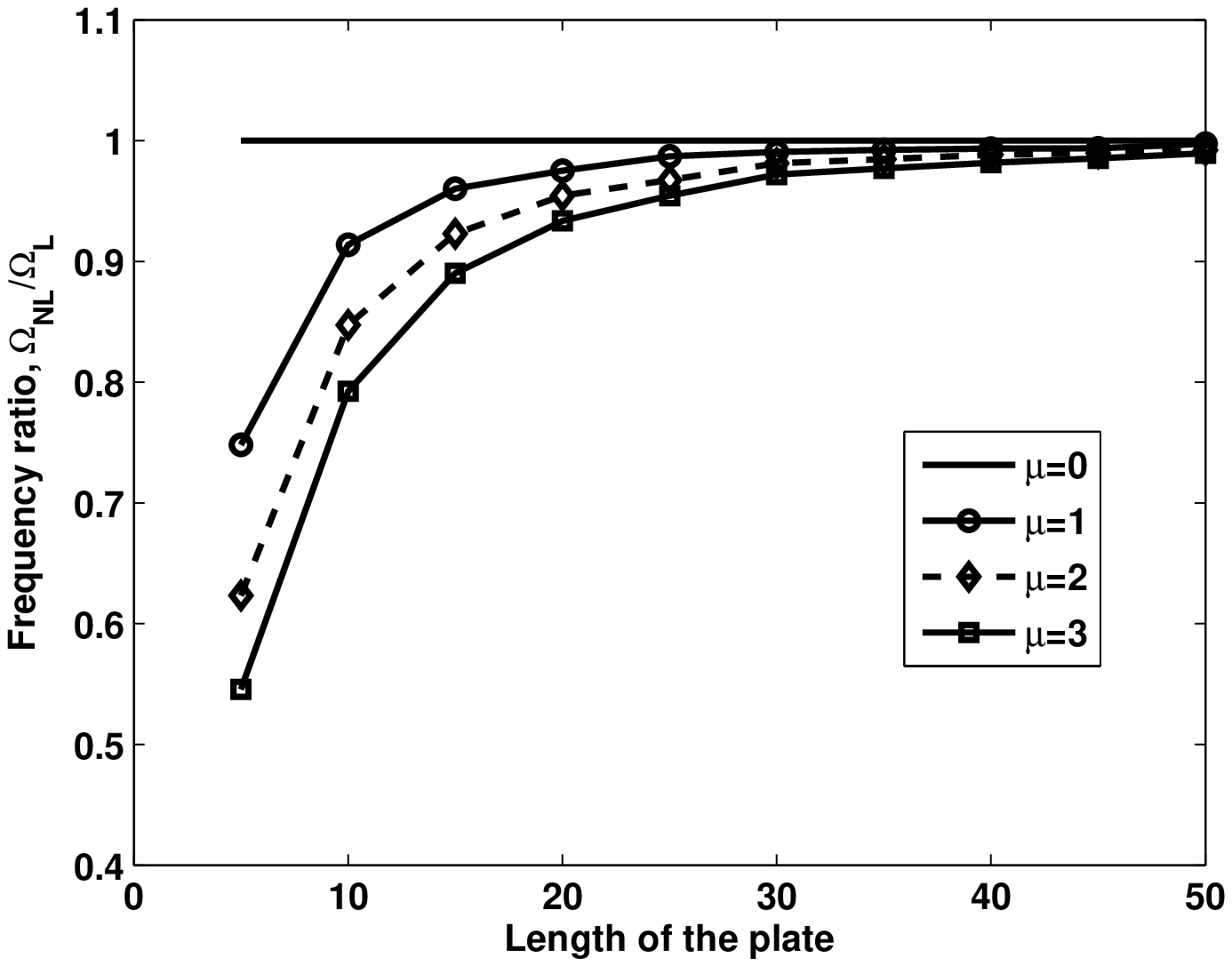}
\caption{Effect of plate dimensions on the frequency ratio $(\Omega_{NL}/\Omega_L)$ for a simply supported square plate for various internal length with $a/h=$ 100.}
\label{fig:plateSSiso_intlen}
\end{figure}

\begin{figure}[htpb]
\centering
\includegraphics[scale=0.6]{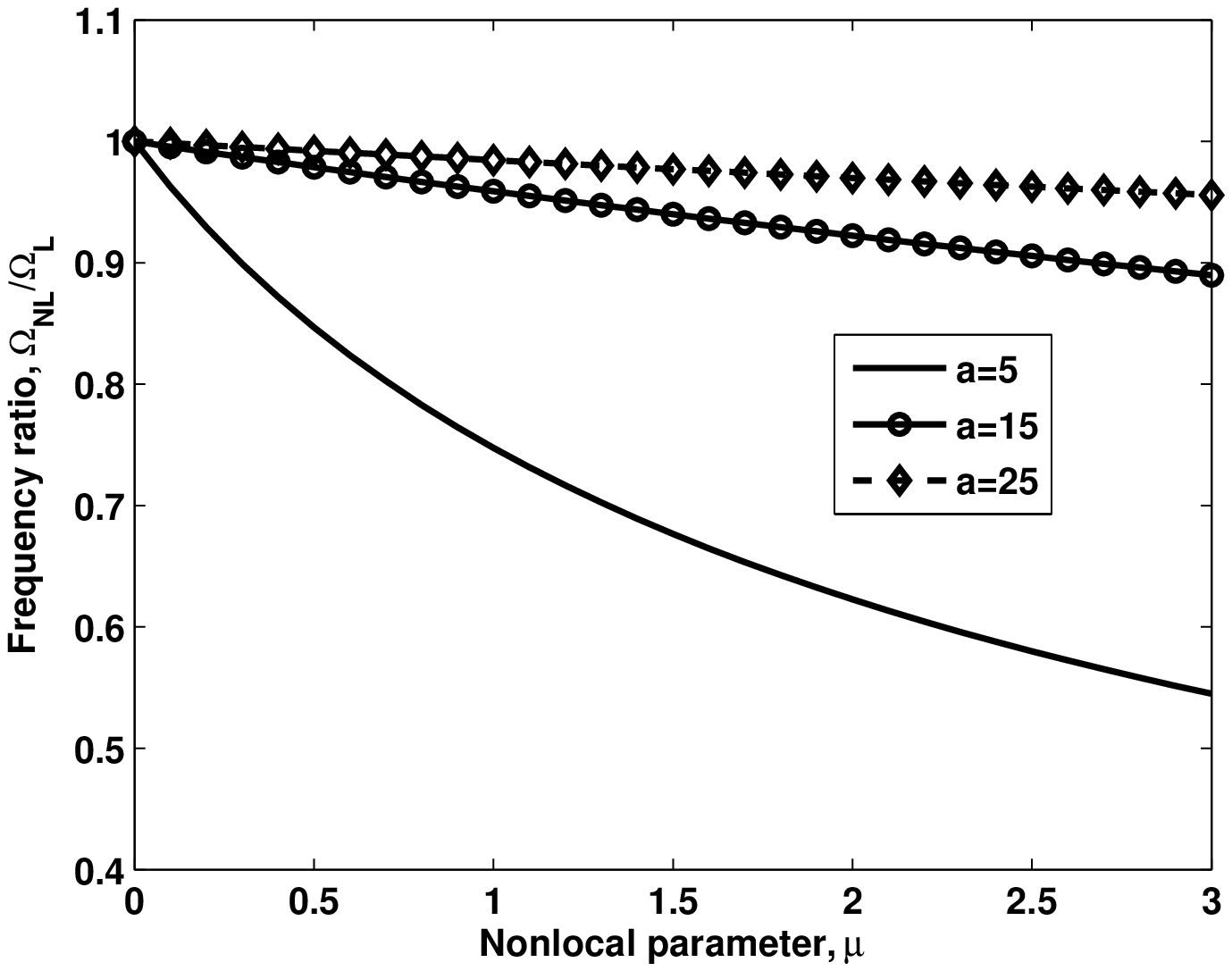}
\caption{Effect of internal length on the frequency ratio $(\Omega_{NL}/\Omega_L)$ for a simply supported square plate for various plate dimensions with $a/h=$ 100.}
\label{fig:plateSSisoLen_intlen}
\end{figure}

\begin{figure}[htpb]
\centering
\includegraphics[scale=0.6]{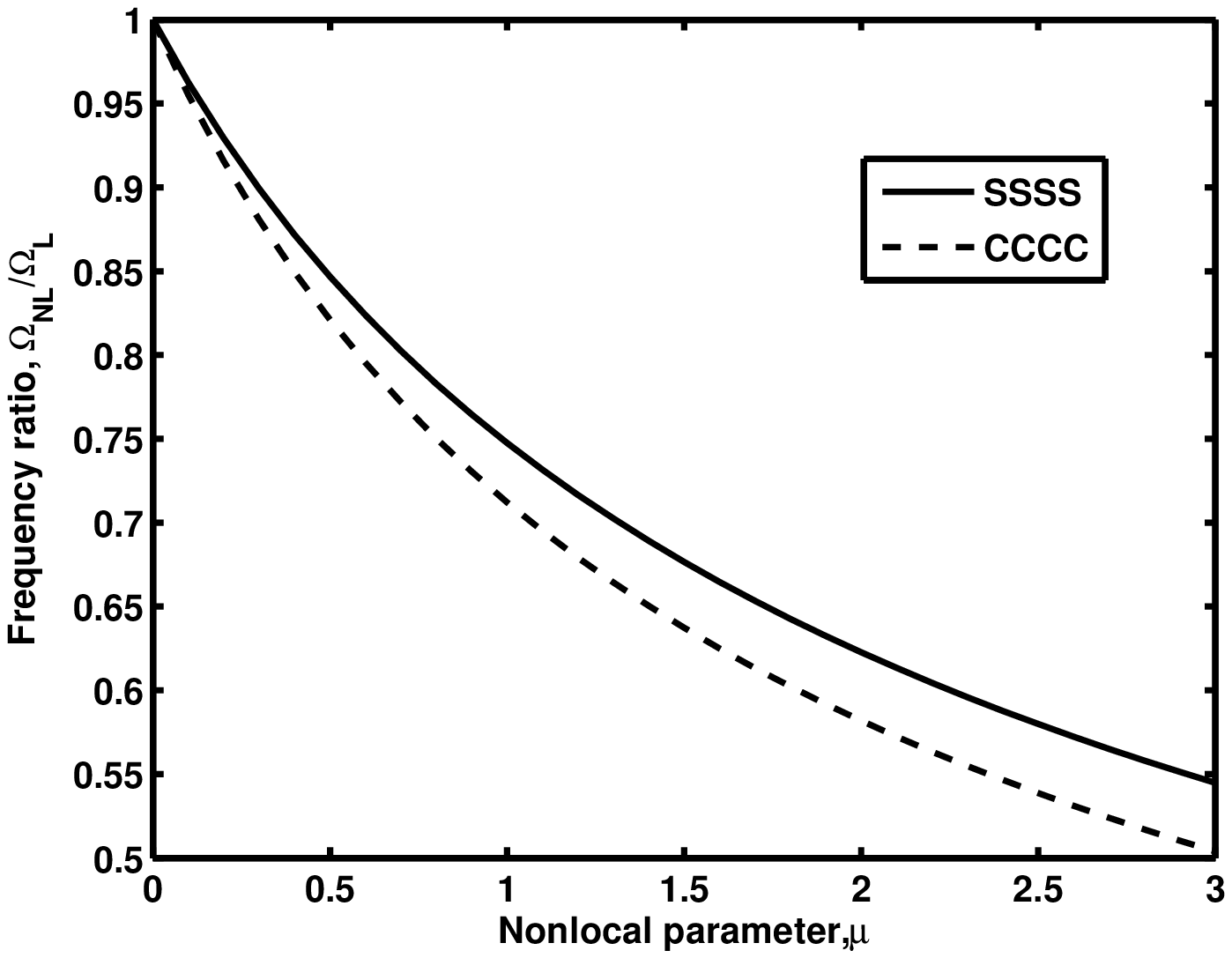}
\caption{Effect of internal length on the frequency ratio $(\Omega_{NL}/\Omega_L)$ for all edges simply supported and all edges clamped boundary condition with $a/h=$ 100 and $a=$ 5.}
\label{fig:boundaryeffec}
\end{figure}

The influence of plate dimension on the frequency ratio for various nonlocal parameter is plotted in \fref{fig:plateSSiso_intlen}. The values of nonlocal parameter are assumed to vary between $\mu=$ 0 (local elasticity) to $\mu=$ 3 nm$^2$. As the length of the plate increases, the frequency ratio tend to increase and approach the local elasticity solution for considerably larger plate length, irrespective of the nonlocal parameter. The influence of nonlocal parameter is significant for smaller plate dimensions. \fref{fig:plateSSisoLen_intlen} shows the influence of nonlocal parameter on the frequency ratio for various plate dimensions. It can be seen that the nonlocal parameter has a stronger influence for smaller plate dimensions, as would be expected. To study the influence of nonlocal parameter and the boundary condition on the natural frequency of an isotropic plate, a square plate of length 10 nm is considered with $a/h=$ 100. The frequency ratio corresponding to two different boundary conditions, viz., all edges simply supported and all edges clamped boundary conditions are plotted in \fref{fig:boundaryeffec}. It can be seen that both the boundary condition and the nonlocal parameter has an influence on the frequency ratio. The higher the nonlocal parameter, the larger is the influence, irrespective of the boundary condition. 


\section{Conclusion}
In this paper, the constitutive model proposed by Eringen is used to model the Timoshenko beam and the plates based on FSDT. The natural frequencies of Timoshenko beam and first order shear deformable plate are studied using Lagrange polynomials and NURBS basis functions. Numerical experiments have been conducted to bring out the effect of boundary condition and the nonlocal parameter on the natural frequency of the beams and plates. The computed results are found to be in good agreement with analytical results. It can be seen that the NURBS basis functions requires fewer degrees of freedom to yield same order of accuracy as that of the RBF in case of the beams. It can be inferred that the effect of nonlocal parameter is to reduce the natural frequency of the plate irrespective of the boundary conditions.

\bibliographystyle{elsarticle-num}
\bibliography{nle_reference}

\end{document}